\documentclass{article}
\usepackage{amssymb}
\usepackage{graphicx}
\usepackage{amsmath}

\newtheorem{theorem}{Theorem}[section]

\newtheorem{corollary}[theorem]{Corollary}

\newtheorem{definition}[theorem]{Definition}

\newtheorem{lemma}[theorem]{Lemma}

\newtheorem{proposition}[theorem]{Proposition}
\newtheorem{remark}[theorem]{Remark}

\newcommand{\SL}{\mathrm{SL}}
\newcommand{\GL}{\mathrm{GL}}
\newcommand{\EL}{\mathrm{EL}}

\newcommand{\mb}[1]{\mbox{\rm {#1}}}
\newcommand{\mr}[1]{\mathcal{#1}}

\newcommand{\R}{\mathbb{R}}

\newcommand{\Z}{\mathbb{Z}}
\newcommand{\N}{\mathbb{N}}
\newcommand{\F}{\mathbb{F}}
\newcommand{\Id}{\mathrm{Id}}
\newcommand{\KC}{\mr{K}}
\newcommand{\TC}{\tau}
\newcommand{\TauC}{$\tau$-constant}
\newcommand{\KaC}{Kazhdan constant}

\begin{document}
\title{Universal lattices and property $\tau$}
\author{M.\ Kassabov and N.\ Nikolov} \date{}

\maketitle
{
\renewcommand{\thefootnote}{}
\footnotetext{
\emph{K}-theory, Steinberg group.}
}
\begin{abstract}
We prove that the universal lattices -- the groups $G=\SL_d(R)$ where $R=\Z[x_1,\dots,x_k]$,
have property $\tau$ for $d\geq 3$. This provides the first example of linear groups with $\tau$
which do not come from arithmetic groups.
We also give a lower bound  for the
\TauC\ with respect to the natural
generating set of $G$.
Our methods are based on bounded elementary generation
of the finite congruence images of $G$, a generalization of a result by
Dennis and Stein on $K_2$ of some finite commutative rings and a
relative property \emph{T} of $(\SL_2(R) \ltimes R^2, R^2)$.
\end{abstract}

{
\renewcommand{\thetheorem}{\arabic{theorem}}
\section*{Introduction}

The groups $\SL_d(\mr{O})$, where $\mr{O}$ is a ring of integers is
in a number field $K$, have many common properties -- they all
have Kazhdan property \emph{T}, a positive solution of the
congruence subgroup problem and super rigidity. In~\cite{YSh},
Y.~Shalom conjectured that many of these properties are inherited from
the group $\SL_d(\Z[x])$. He called the groups
$\SL_d(\Z[x_1,\dots,x_k])$ \emph{universal lattices}, because they
can be mapped onto many lattices in groups $\SL_d(K)$ for different
fields $K$. Almost nothing is
known about the representation theory of these groups.

The main result in this paper is Theorem~\ref{main}, which says that
the universal lattices have property $\tau$, provided that $d\geq 3$.
However unlike the classical lattices $\SL_d(\mr{O})$ these groups
do not have congruence subgroup property and
have infinite (even infinitely generated) congruence kernel.

\medskip

Let $G$ be a topological group and consider the space
$\widetilde{G}$ of all equivalence classes of unitary
representations of $G$ on some Hilbert space $\mr{H}$. This space
has a naturally defined topology, called the Fell topology, as
explained in \cite{LubZuk}  \S1.1 or \cite{expanderbook} Chapter 3 for example.
Let $1_G$ denote the
trivial $1$-dimensional representation of $G$ and let
$\widetilde{G}_0$ be the set of representations in $\widetilde{G}$
which do not contain $1_G$ as a subrepresentation
(i.e., do not have invariant vectors).

\begin{definition}
A group is $G$ is said to have Kazhdan property
\emph{T} if $1_G$ is isolated from
$\widetilde{G}_0$ in the Fell topology of $\widetilde{G}$.
\end{definition}

A discrete group with property \emph{T} is finitely generated: see
\cite{kazhdan}. Since in this paper we shall be concerned only
with discrete groups we give the following equivalent
reformulation of property \emph{T}:

\emph{Equivalent definition:}
Let $G$ be a discrete group generated by a finite set $S$.
Let $\rho: G \rightarrow U(\cal{H})$ be a unitary representation
of the group $G$. A vector $v\in \mr{H}$ is called $S$
$\epsilon$-invariant vector iff $||\rho(s)v-v||< \epsilon ||v||$ for all
$s\in S$. Then $G$ has the Kazhdan property \emph{T} if there is
$\epsilon  >0$ such that every irreducible  unitary
representation $\rho: G \rightarrow U(\cal{H})$ on a
Hilbert space $\mr{H}$, which contains a $\epsilon$-invariant vector
is isomorphic to the trivial representation.
The largest $\epsilon$ with this property is called the \emph{Kazhdan Constant}
for $S$ and is denoted by $\KC(G;S)$.
\medskip

The property \emph{T} depends only on the group $G$
and does not depend on the choice of the generating set $S$, however the
Kazhdan constant depends also on the generating set.

Property \emph{T} implies certain group theoretic conditions on $G$
(finite generation, FP, FAB etc) and can be used for construction
of \emph{expanders} from the finite images of $G$.
For this last application the following weaker property $\tau$ (introduced
by Lubotzky in~\cite{expanderbook}) is sufficient:

Let $\widetilde{G}^f$ and $\widetilde{G}^f_0$ denote the finite representations
of $\widetilde{G}$, resp. $\widetilde{G}_0$ (i.e., the representations which factor
through a finite index subgroup).

\begin{definition}
A group is $G$ is said to have property $\tau$ if $1_G$
is isolated from $\widetilde{G}^f_0$
in the induced Fell topology of $\widetilde{G}^f$.
Equivalently: the group $G$ with the
profinite topology
has property \emph{T}.
\end{definition}
The two definitions are equivalent because any continuous irreducible
representation of $G$ in $U(\mr{H})$ is then finite.

Again, for a discrete finitely generated group $G$ there
is an equivalent definition of property $\tau$:
\medskip

Let $G$ be an discrete group generated by a finite%
\footnote{There are examples of groups with property $\tau$ which
are not finitely generated. This definition can be modified to
cover also this case. We are not going to need this because by
a result of Suslin~\cite{suslin} the universal lattices are finitely
generated.}
set  $S$.
Then $G$ has the property $\tau$ if there exists
$\epsilon  >0$ such that for every finite irreducible  unitary
representation $\rho: G \rightarrow U(\cal{H})$ on a
Hilbert space $\mr{H}$ (necessary finite dimensional),
which contains an $\epsilon$-invariant vector
is trivial.
The largest $\epsilon$ with this property is called the \emph{\TauC}
and is denoted by $\TC(G;S)$.

\medskip

Our approach to property $\tau$ is inspired by a paper by
Shalom \cite{YSh} which relates property \emph{T} to
\emph{bounded generation}.
In this paper we will work only with the groups
$\SL_d(R)$, where $d>2$ and $R$ is a finitely generated
commutative ring. The arguments
can be easily generalized to any high rank Chevalley group
over $R$. It is also possible to extend some parts
of the argument to Chevalley groups over noncommutative
rings~\cite{KSL3k}.

\medskip

Let $R$ be a commutative (unital) ring, and for $i \not = j \in \{1,2,\dots,d\}$ let
$E_{i,j}$ denote the set of elementary $d \times d$ matrices
$\{\Id+r\cdot e_{i,j}\  | \ r \in R \}$. Also set
$E=E(R)=\bigcup_{i\not =j}E_{i,j}$ and
let $\EL(d;R)$ be the subgroup in $\GL_n(R)$ generated by $E(R)$.
By a result of Suslin~\cite{suslin}  we have
that $\SL_d(R)=\EL(d;R)$, i.e., that $SK_1(d,R)=1$ in the case of
$d\geq 3$ and $R=\Z[x_1,\dots,x_k]$;
for $d=2$ this is not true in general.

\begin{definition}
The group $G=\EL(d;R)$ is said to have
bounded elementary generation property
if there is a
number $N=BE_d(R)$ such that every element of $G$
can be written as a product of at most $N$
elements from $E(R)$.
\end{definition}

Examples of $R$ satisfying the above definition are rings of
integers $\mr O$ in number fields $K$ (for $d\geq 3$), see \cite{CK}.
In this classical case this property is
known as \emph{bounded generation} because each group $E_{i,j}\simeq (\mr O,+)$
is a product of finitely many cyclic groups. \medskip

The following theorem was proved in~\cite{YSh} and the
method of its proof forms the basis of our results .
\begin{theorem}
\label{t1}
Suppose that $d\geq 3$, $R$ is a $k$-generated commutative ring
such that $\SL_d(R)=\EL(d;R)$ has bounded elementary generation property.
Then $\SL_d(R)$ has property \emph{T} (as a discrete group).
Moreover the Kazhdan constant $\mr{K} (G,S)$ is bounded from bellow by
$$
\mr{K} (G,S) \geq \frac{1}{BE_d(R) 22^{k+1}},
$$
for a specific generating set $S$ (defined below).
\end{theorem}
\medskip

The generating set $S$ in Theorem \ref{t1} is defined as follows:
Suppose that the ring $R$ is generated by $1$ and $\alpha_{1},\dots,\alpha_{k}\in R$. Then
$S=S_{d,k}:=F_{1}\cup F_{2}$, where
$F_{1} =\{\Id \pm e_{i,j}\}$ is the set of $2(d^2-d)$ unit elementary matrices, and
$$
F_{2}=\left\{\Id\pm\alpha_{l}\cdot e_{i,j} \ | \ |i-j|=1, 1\leq l\leq k\right\},
$$
is the set of $4(d-1)k$ elementary matrices with generators of the ring $R$
next to the main diagonal.

\bigskip

A very interesting conjecture is whether the group
$G_{d,1}:=\SL_d(\Z[x])$ and the other universal lattices
have bounded elementary generation property.
In view of Theorem \ref{t1} this would imply that
$G_{d,1}$ (and therefore all of its images which include $\SL_d(\mr{O})$
for many rings of algebraic integers $\mr{O}$%
\footnote{Property \emph{T} for each of the groups $\SL_d(\mr{O})$ is known, but
the Kazhdan constant depends slightly on the ring $\mr{O}$. Property \emph{T} of
$\SL_d(\Z[x])$ would give a uniform property \emph{T} for almost all $\SL_d(\mr{O})$.
}) %
 has property \emph{T}.
We are unable to say anything about property \emph{T},
but  we shall prove that $G_{d,1}$ at least has $\tau$.
More generally:
\begin{theorem}
\label{main}
Let $d\geq 3$, $k \geq 0$ denote the ring
$\Z[x_1,\ldots, x_k]$ by $R_k$.
Then the universal lattice $G_{d,k}:=\SL_d(R_k)$ has property $\tau$.

Moreover we have the following explicit bound for the
\TauC\ $\TC(G_{d,k};S)$ with respect to the
generating set $S_{d,k}$ as above:%
\footnote{The bound for the Kazhdan constant can be significantly improved, see~\cite{KSL3k}.}
$$
\TC(G_{d,k};S) > \frac{1}{\left(2d^2 + 18k + 30\right)22^{k+1}}
$$
\end{theorem}

It is a nontrivial fact proved by Suslin~\cite{suslin} that the group $G_{d,k}$
is generated by the set $S_{d,k}$. However this is independent from our result because
a group $G$ can have property $\tau$ and a positive \TauC\
with respect to some set $S$ which does not generate the group, it is sufficient
that $S$ generates a dense subgroup (in the profinite topology) in $G$.

Let $\widehat{G_{d,k}}$ be the profinite completion of the universal lattice $G_{d,k}$ and define the
\emph{pro-elementary subgroup} $\widehat{E_{i,j}}$ to be the closure of $E_{i,j}$ in $\widehat{G_{d,k}}$.
Each $\widehat{E_{i,j}}$ is isomorphic to the additive group of the
profinite completion $\widehat{R}$ of the ring $R$. \medskip

We shall prove 
\begin{theorem}
\label{BpEG}
The profinite completion $\widehat{G_{d,k}}$ is `boundedly pro-elementary' generated: it is a finite product
of the groups $\widehat{E_{i,j}}$.

In fact $\widehat{G}$ can be written as a product of at most
$(3d^2-d-2)/2 + 18(k+2)$ pro-elementary subgroups $\widehat E_{i,j}$ 
in some fixed order.
\end{theorem}

In order to prove this Theorem we need Lemma~\ref{congimage}
below on bounded elementary
generation of $\SL_d$ over a finite ring and the following result
which may be of independent interest:

\begin{theorem}
\label{k2}
Let $\bar{R}$ be a finite commutative ring generated by $k$ elements.
Then every element of $K_2(\bar{R})$
is a product of at most $k+2$ Steinberg symbols $\{a,b\}$.
\end{theorem}

Once Theorem \ref{BpEG} has been proved, the general techniques from~\cite{YSh} are
applied to prove that $\widehat{G}$ has property \emph{T}.
As noted above this gives that $G$ has property $\tau$.
\medskip

Most of the known examples of discrete finitely
generated groups $G$ with property $\tau$
arise as lattices in higher rank semi-simple Lie groups.
In particular they almost have property T and `rigidity' (even super-rigidity):
their representation theory is controlled by the representations
of the ambient Lie group. There are also many examples of
`randomly presented' hyperbolic
groups with property T, e.g.\ \cite{gromov,zuk}, but it is not known
if their profinite completion is infinite.

The universal lattices seem to be the first residually finite
`non-arithmetic' groups%
\footnote{In some sense these groups are arithmetic subgroups of $\SL_d(K)$
for some huge field $K$ like $\R((x_1^{-1}))\dots((x_k^{-1}))$,
however they are not cocompact subgroups and
the standard theory of arithmetic groups can not be applied in this case.}
with property $\tau$ discovered so far. They have infinitely (continuously)
many irreducible (but not unitary) representations of fixed finite degree.
In this light the question whether they have property \emph{T} is even more interesting.
\bigskip

Finally we remark that Theorem~\ref{main} gives a lower bound for the \TauC\
which is asymptotically $O(d^{-2}22^{-k})$ in $d$ and $k$.
However it is possible to improve this estimate to
$O\left(d^{-1/2}(1+(k/d)^{3/2})^{-1}\right)$.
The proof of this result involves some new ideas from \cite{K},
namely relative property T of of groups with big radical and use of
\emph{generalized elementary matrices}.
If the reader is interested in the values of the
Kazhdan constants and generalization of Theorem~\ref{main} to non commutative rings,
he is encourage to go over the sequel of this paper~\cite{KSL3k}.
It is interesting to note that using a different generating set it is
possible to improve the the \TauC\ for $G_{d,k}$ to $O((1+k/d)^{-1/2})$.

\bigskip

A few words about the structure of the rest of the paper:

Section~\ref{outline} contains the proof of Theorem~\ref{main} modulo Lemma \ref{congimage}
(proved in Section~\ref{reduction}), Theorem~\ref{k2} (proved in Section~\ref{symbols})
together with some technical results from~\cite{YSh}.
We conclude with Section~\ref{CSP} where we investigate the implication Theorem~\ref{k2} has to
the structure of the congruence kernel of $G$ and prove Theorem~\ref{BpEG}.
\bigskip

\textbf{Acknowledgements: }
We wish to thank R.~K.~Dennis, A.~Lubotzky, Y.~Shalom and L.~Vaserstain for
useful discussions and for the inspiration from their results.

Part of this work was completed while the first author was at the postdoctoral fellow at
University of Alberta, Edmonton and the second one held a visiting research fellowship at TIFR, Mumbai.
\bigskip

\textbf{Notation:} For the rest of this paper we will denote
$R_k:=\Z[x_1,\dots,x_k]$ and $G_{d,k}:=\SL_d(R_k)$.
We will fix a generating set $S_{d,k}=F_1\cup F_2$ of the group $G_{g,k}$
where
$F_1$ is the set of $2(d^2-d)$ elementary matrices with $\pm 1$ off the diagonal
and $F_2$ is the set $4k(d-1)$ elementary matrices $\Id \pm x_ne_{ij}$ with $|i-j|=1$.
Unless we explicitly need to specify the number of generators of the ring
$R_k$ or the size of matrices $G_{d,k}$
we will denote $R_k$ by $R$, $G_{d,k}$ by $G$ and $S_{d,k}$ by $S$.
We will assume that all rings are commutative with $1$.
}

\section{Proof of Theorem~\ref{main} }
\label{outline}

Let $\rho : G_{d,k} \to U(\mr{H})$ be a unitary representation of the
universal lattice $G_{d,k}$ which factors through a
finite index subgroup. Suppose that $v\in \mr{H}$ is an $\epsilon$-almost
invariant unit vector for the set $S$, where
$\epsilon \leq  K_{d,k}$

The following result is Corollary 3.5 to Theorem 3.4 proved in~\cite{YSh}.
\begin{lemma}
\label{abe}
There is a constant $M(k)<11.22^{k}$ with the following property:
Let $(\rho ,\mr{H})$ be a unitary representation of the group
$G_{d,k}=\SL_d(R_k)$.
Let $v\in \mr{H}$ be a unit vector such that $||\rho(g)v-v|| <\epsilon$ for $g\in S_{d,k}$.
Then for every elementary matrix $g = \Id + r e_{ij}$
we have $||\rho(g)v - v || < 2M(k)\epsilon$.
\end{lemma}

We shall prove that the \TauC\ for the universal lattice $G_{d,k}$
is bounded from below by:
$$
\TC(G_{d,k};S_{d,k}) \geq K_{d,k}:=[(3d^2-d-2) + 36(k+2)]^{-1}M(k)^{-1}.
$$
\medskip

Applying Lemma~\ref{abe} gives us that
$$
||\rho(g)v-v|| < 2M(k) \epsilon
$$
for every elementary matrix $g\in E$.

\bigskip

Let $H = \ker \rho < G$, this is a normal subgroup of finite index in $G$.
Define subgroups $H_{i,j}$ of $H$ and subsets $U_{i,j}$ of the ring $R$ for $i \not =j$ as follows
$$
H_{i,j}:= H \cap E_{i,j}; \quad \quad U_{i,j}:=\{r\in R | \,\,\, 1 + r\cdot e_{i,j} \in H\}.
$$
Using elements in the Weyl group we see that the subgroups $E_{i,j}$ are pairwise conjugate in $\SL_d(R)$.
This gives us that the sets $U_{i,j}$ do not depend on the indices $i,j$.
Using commutation with a suitable elementary matrix we can see that $U_{i,j}$ is also closed under
multiplication by elements in $R$ (this is valid only for $d\geq 3$), i.e., $U=U_{i,j}$ is an ideal
of the ring $R$. This ideal is of finite index in $R$ because $H$ is a subgroup of finite index in $G$.

Let $\EL(d;U)$ be the normal subgroup in $G$ generated by $H_{i,j}$:
$$
\EL(d;U):= \langle H_{i,j} \rangle^G < H
$$

Let $SK_1(R,U;d)=G(U)/\EL(d;U)$ where $G(U)$ is the principal congruence subgroup of
$G$ modulo $U$, i.e., the matrices in $G$ congruent to $\Id_d$ modulo $U$.
We have the following diagram
$$
\begin{array}{c@{}c@{}c@{}c@{}c@{}c@{}c@{}c@{}c}
1 & \rightarrow & SK_1(R,U;d) & \rightarrow & G/\EL(d;U) & \rightarrow & G(R/U)= \SL_d(\bar{R}) & \rightarrow & 1 \\
  &             &             &             & \downarrow&             &                        &             & \\
  &             &             &             &   G/H     &             &                        &             & \\
  &             &             &             & \downarrow&             &                        &             & \\
  &             &             &             &    1      &             &                        &             & \\
\end{array}
$$
where the row and the column are exact.

The next lemma is well known in the case of local fields, for general commutative rings is
slightly stronger than Remark~10 in~\cite{DV}:
\begin{lemma}
\label{congimage}
Let $\bar R$ be a finite commutative ring.
Then $\SL_d(\bar{R})$ has uniform bounded elementary generation property.
In fact every matrix in
$\SL_d(\bar{R})$ can be written as a product of $(3d^2-d-2)/2$
elementary matrices.
\end{lemma}

It remains to deal with the group $SK_1(R,U;d)$. It is always a finite abelian group and
it measures the departure of $\EL(d;U)$ from being
a congruence subgroup. We note that while $SK_1(R,U;d)$ is finite,
it is not bounded in terms of $d$ and $R$ and can be
arbitrarily large for some ideals $U$ (unless $R=\Z$, when $SK_1(R,U;d)$ is always trivial),
see Section~\ref{CSP}.
\medskip

Nevertheless we have the following:
\begin{lemma}
\label{kernel}
Each element of the group $G(U)/\EL(d;U)=SK_1(R,U;d)$ is a product
of at most $18(k+2)$ elementary matrices.
In fact it is a product of at most $k+2$ Steinberg symbols $\{a,b\}$,
for $a,b$ invertible elements in $\bar{R}$,
coming from the surjection
$K_2(\bar R) \rightarrow SK_1(R,U;d)$ (to be defined in Section~\ref{symbols}).
\end{lemma}
We obtain this Lemma as a corollary to Theorem~\ref{k2} proved in Section~\ref{symbols}.

To finish the proof of Theorem~\ref{main} we use an argument due to
 Y.~Shalom from~\cite{YSh}:
 Let $g$ be any element in the group $G$. Using Lemmas~\ref{abe} and~\ref{kernel} every
$g\in G$ can be written as a product
$$
g = h\prod_{s=1}^{N} u_s
$$
where $h\in \EL(d;U) \subset H$ and $u_s\in E(R)$ are elementary matrices, here
$N$ denotes the number $(3d^2-d-2)/2 + 18(k+2)$.
By Lemma \ref{abe}
$$
\begin{array}{c@{}c@{}l}
\displaystyle
||\rho(g)v-v|| & {}\leq {}& \displaystyle ||\rho(h)v-v|| + \sum_{s=1}^{N} ||\rho(u_s)v-v|| < \\
               & < & \displaystyle 0 + \sum_{s=1}^{N} 2M(k)\epsilon  = 
               \frac{(3d^2-d-2) + 36(k+2)}{2} 2 M(k) \epsilon.
\end{array}
$$

Recall that $\epsilon \leq  K_{d,k}$ i.e.,
$$
||\rho(g)v-v|| < \frac{(3d^2-d-2) + 36(k+2)}{2} 2 M(k) K_{d,k} = 1.
$$
This shows that every element in the group $G$ moves the vector $v$ by less then $1$.
Let $V$ denote the the $G$-orbit of the vector $v$ in $\mr{H}$. The center of mass $c_V$
of the set $V$ is invariant under the action of $G$ and is not zero because the whole orbit $V$
lies entirely in the half space $\{v'\in \mr{H} \mid \Re (v,v') > 0\}$.
Therefore $\mr{H}$ contains a nontrivial $G$-invariant vector, which completes the
proof that the group $G$ has property $\tau$.

Finally, the estimate for $M(k)$ show that the constant $K_{d,k}$ satisfies the inequality
$$
K_{d,k} \geq \frac{1}{\left(2d^2 + 18k + 30\right)22^{k+1}}.
$$
Theorem~\ref{main} is now proved, modulo Theorem~\ref{k2} and Lemma~\ref{congimage}.

\section{Uniform bounded generation of $\SL_d (\bar{R})$}
\label{reduction}

Lemma \ref{congimage} is a well-known result for fields
and the same proof works for any local ring $\bar R$:
Where the classical algorithm
(described below)
uses a nonzero `pivot' field element in a row or column operation, the algorithm for
$\bar R$ uses the corresponding element outside the maximal ideal $\bar I$ of $\bar R$,
which is therefore a unit.
This argument can be generalized to semi-local rings if one is careful about the order of
the elementary matrices which appear in the product. See also \cite{DV}, where a more
general result is proved for rings with Bass stable range at most 1.

It is well known that any
finite commutative ring is semi-local. In fact:

\begin{theorem}
\label{rings}
Any commutative unital finite ring is a direct sum of local rings.
\end{theorem}
\textbf{Proof:}
This is known as a generalized Chinese Reminder Theorem. The proof is by induction on
the size of the finite ring $\bar R$: Let $\bar J$ be the Jacobson radical $\bar R$.
Then $\bar R/\bar J$ is a finite direct product of finite fields.
If it is a field there is nothing to prove because $\bar J$ is a maximal ideal and $\bar R$
is a local ring.
Otherwise let $e_0\in \bar R$ map to a nontrivial
idempotent in $\bar R/\bar J$.
Since $\bar J$ is nilpotent ideal we can find $e \in \bar R$ such that
$e \equiv e_0(\mb{mod } \bar J)$ and $e$ is an
(nontrivial) idempotent in  $\bar R$. Then we have
$\bar R = e\bar R \oplus (1-e)\bar R$.
Finally we can apply the induction hypothesis to the rings $e\bar R$ and $(1-e)\bar R$.
$\square$
\medskip

It is clear that if $\bar R$ is a direct sum
$\bar{R}=\bar{R}_1 \oplus \bar{R}_2 \oplus \cdots \oplus \bar{R}_s$
then the linear group $\SL_d(\bar R)$ decomposes as a direct product
$$
\SL_d(\bar{R}) \simeq \SL_d(\bar{R}_1) \times \cdots \times \SL_d(\bar{R}_s).
$$
This observation together with Theorem \ref{rings} reduces the proof of Lemma \ref{congimage}
to the case of local rings $\bar{R}$, namely the following
Lemma~\ref{congimagest}:

\begin{lemma}
\label{congimagest}
Let $\bar R$ be a local ring.
Then $\bar{G}=\SL_d(\bar{R})$ has bounded generation with respect to
$\bigcup_{i,j} E_{i,j}(\bar{R})$. In fact $\bar G$
can be written as
$$
\bar G= \prod_{s=1}^{(3d^2-d-2)/2}E(s),
$$
where each of the $(3d^2-d-2)/2$ terms $E(s)$ in the product is
one from the groups $E_{i,j}(\bar R)$ (with entries in $\bar{R}$) in
some fixed order.
\end{lemma}
\textbf{Proof:}
Since $\bar{R}$ is a local ring Lemma \ref{congimagest} is a consequence of the familiar argument
that a matrix $g \in \SL_d(K)$ can be reduced by successive applications of row and column
operation to the identity matrix: Each of these operations is in fact a multiplication by
an elementary matrix from left or right.
\medskip

The following well-known algorithm produces such decomposition%
:

Let $g\in \SL_d(\bar R)$. By
multiplying with $d-1$ elementary matrices to the right we can ensure that the last entry
on the first row is an invertible element in $\bar R$. Using an extra multiplication to the
right we can make the $d,d$ entry  equal to $1$, next with $d-1$ left and $d-1$ right multiplications
by elementary matrices we can transform the matrix $g$ to an element in
$\SL_{d-1}\times 1$ (sitting in the top left corner of $\SL_d$) .
Thus the reduction form $\SL_d$ to $\SL_{d-1}$ can be done using $3d-2$ elementary matrices,
and by induction using $(3d^2-d-2)/2$ we can transform any matrix to the identity.

This proves Lemma \ref{congimagest} and Lemma \ref{congimage} follows.
$\square$



\section{$K_2$ of finite commutative rings: Theorem \ref{k2}}
\label{symbols}
In this section we allow $\bar{R}$ to be any commutative ring with 1. Let $d>2$.
In order to simplify the argument a bit we require that $\mathrm{SL}_d(\bar{R})$ is
generated by the elementary matrices $E_{i,j}(r)$ with entries in $\bar{R}$.
This assumption is indeed true for any finite ring $\bar R$
(Lemma~\ref{congimage} is even stronger).
Recall the definition of the Steinberg group $\mathrm{St}_d(\bar{R})$:%
\footnote{Analogues of the Steinberg group can be defined for the
other Chevalley groups, see~\cite{stein} and the references therein.
The results from this section
also hold in that case.}

\begin{definition}
$\mathrm{St}_d(\bar{R})$ is the group generated by elements $x_{i,j}(r)$ for
distinct $i,j \in \{1,2,\dots,d\}$ and $r \in \bar{R}$ subject to the relations:
\begin{enumerate}
\item $x_{i,j}(u)x_{i,j}(v)=x_{i,j}(u+v)$,
\item $[x_{i,j}(u),x_{j,k}(v)]=x_{i,k}(uv)$,
\item $[x_{i,j}(u),x_{k,l}(v)]=1$
\end{enumerate}
for all distinct indices $i,j,k,l \in \{1,2,\dots,d\}$ and all $u,v \in \bar{R}$.
\end{definition}

We have an obvious surjection
$$
\phi: \ \mathrm{St}_d(\bar{R}) \rightarrow \EL(d;\bar{R})=\SL_d(\bar{R})
$$
given by $x_{i,j}(u)\mapsto \Id  + u e_{i,j}(u)$.
The kernel of $\phi$ is of great interest in classical $K$-theory and related areas and is denoted by
$K_2(\bar R;d)$.%
\footnote{For a finite ring $\bar R$ it is in fact independent of $d\geq 3$ by stability
results in \cite{bass}, but we shall not need this.}
Unless we need to specify the size of the matrices
we will denote $K_2(\bar R;d)$ by $K_2(\bar R)$.
\bigskip

Some naturally defined elements of it are defined below:

For a unit $u$ in $R$ and a pair of indices $i,j$ consider the elements
$$
w_{i,j}(u):=x_{i,j}(u)x_{j,i}(-u^{-1})x_{i,j}(u) \textrm{ and } h_{i,j}(u):=w_{i,j}(u)w_{i,j}(-1).
$$
For any units $u,v$ we define the Steinberg symbol
$$
\{u,v\}_{i,j}=h_{i,j}(uv) h_{i,j}(u)^{-1}h_{i,j}(v)^{-1} \in \mathrm{St}_d(\bar{R}).
$$
It can be shown that the element $\{u,v\}_{i,j} \in \mathrm{St}_d(\bar{R}) $ is independent on the choice of
$i,j$ and lies in the kernel of $\phi$. This element in $K_2(\bar{R})$ is called a
Steinberg symbol and is denoted by $\{u,v\}$. It has the properties 
\begin{enumerate}
\item bimultiplicative: $\{uv,w\}=\{u,w\}\{v,w\}$,  
\item skew-symmetric: $\{u,v\}=\{v,u\}^{-1}$,       
\item $\{u,-u \}=1$ and                             
\item $\{u,1-u\}=1$                                 
\end{enumerate}
for all units $u,v,w$ (the last identity only if $1-u$ is also a unit).
See~\cite{milnor}, Chapter 9 for details.

\begin{theorem}
[Dennis and Stein, \cite{DS2}]
When $\bar{R}$ is
a semi-local commutative ring then $\ker \phi=K_2(\bar{R};d)$ is
central in $\mathrm{St}_d(\bar{R})$ and is generated by the
Steinberg symbols $\{a,b\}$ for all $a,b$ invertible in $R$.
\end{theorem}

Let $\pi: \SL_d(R)/\EL(d;U) \rightarrow \SL_d(\bar{R})$ be the reduction modulo $U$.
Then $\phi$ factors through $\pi$.

Indeed, define a map
$\psi: \mathrm{St}_d(\bar{R}) \rightarrow \SL_d(R)/\EL(d;U)$ by
$x_{i,j}(\bar{r}) \mapsto \Id+r \cdot e_{i,j}$.
This is well defined and extends to a homomorphism $\psi$ with the property that $\phi=\pi \circ \psi$.
We thus see that $\ker \pi =\SL_d(U)/\EL(d;U)=SK_1(R,U;d)$ is an image of $K_2(\bar{R};d)$ under $\psi$
and is therefore central in $\SL_d(R)/\EL(d;U)$ and generated by the images $\psi (\{a,b\})$ of symbols.
Another way to see this is to use the exact sequence from classical
\emph{K}-theory cf. Theorem 6.2 in \cite{milnor}:
$$
\cdots \rightarrow K_2(R;d) \rightarrow K_2(\bar R;d) \stackrel{\psi}{\rightarrow} SK_1(R,U;d)
\rightarrow SK_1(R;d) \rightarrow SK_1(\bar R;d) \rightarrow \cdots
$$
and by the result of Suslin \cite{suslin} $SK_1(R) = SK_1(\bar R)=1$, hence our map $\psi$ is surjective.
\medskip

By definition each $\psi (\{a,b\})$ is a product of 18 images of
elementary matrices in $\SL_d(R)/\EL(d;U)$ (it is actually a product of only $13$ images of
elementary matrices because there are some cancelations)
and therefore Lemma \ref{kernel} follows from Theorem \ref{k2} in the Introduction.

The rest of this section is devoted to the proof of Theorem \ref{k2}.
\bigskip

First, note that the functor $K_2$ respects direct sums of rings i.e.
$$
K_2(R_1 \oplus \cdots \oplus R_s)=K_2(R_1)\oplus \cdots \oplus K_2(R_s).
$$
Hence by Theorem \ref{rings} we may assume that $\bar{R}$ is a
finite local ring with a maximal ideal $\bar I$.

\begin{lemma}
If the finite local ring $\bar R$ is an image of
$R_k=\Z[x_1\,\dots,x_k]$ then the maximal ideal $\bar I$ in $\bar{R}$
is generated by at most $k+1$ elements.
\end{lemma}
\textbf{Proof:}
Indeed $\bar I$ is an image of a maximal ideal $I$ of finite index in $R_k=\Z[x_1,\dots,x_k]$.
If $p$ is the characteristic of $R/I$ then $I/pR$ is a maximal ideal of
$\F_p[x_1,\dots,x_k]$ and Theorem~24 of Chapter~VII in~\cite{ZS} gives that $I/pR$ is a
$k$-generated ideal. This proves the claim. There is another, more
conceptual, proof of this lemma based on the fact that any finite field has a presentation
as a ring with $2$ generators and $2$ relations.
$\square$
\medskip

In their paper \cite{DS1} Dennis and Stein proved several new identities
for the Steinberg symbols besides the standard one we quoted:

\begin{theorem}
[\cite{DS1} Proposition 1.1]
\label{identities}
Let $v,p,q,r$ be elements of $\bar R$ such that
$v, 1-p,1-q,1-r, 1-qv, 1-pv, 1-pqv, 1-pq, 1-pr, 1-qr, 1-pqr$ are units in $\bar R$.
Then the following identities hold in $K_2(\bar R)$:
\begin{equation}
\label{I1}
\{v,1-pqv\}=
\left\{  -\frac{1-qv}{1-p},  \frac{1-pqv}{1-p} \right\}
\left\{  -\frac{1-pv}{1-q},  \frac{1-pqv}{1-q} \right\}
\end{equation}
and
\begin{equation}
\label{I2}
 \left\{  -\frac{1-qr}{1-p},  \frac{1-pqr}{1-p} \right\}
 \left\{  -\frac{1-pr}{1-q},  \frac{1-pqr}{1-q} \right\}
 \left\{  -\frac{1-pq}{1-r},  \frac{1-pqr}{1-r} \right\}
=1
\end{equation}
\end{theorem}
\medskip

The basic idea of the following construction is derived from~\cite{DS1}, where the authors
prove the $K_2(\bar R)$ is a cyclic group if $\bar R$ is a factor of a discrete valuation
ring. We can not apply directly this result because we need a result which holds for any finite
local ring and there are many finite local
rings which do not arise from discrete valuation rings.
\medskip

For a nonzero element $a\in \bar R$ define the \emph{level} $l(a)$ to be the
largest $n$ such that $a \in \bar{I}^n$ and set $l(0) =\infty$.
In particular $l(a)=0$ if and only $a$ is a unit of $\bar R$.

For $n\in \N$, let $K_2^{(n)}$ be the subgroup of $K_2(\bar{R};d)$
generated by all symbols $\{a,b\}$ with $a,b$ units such that
$l(a-1)+l(b-1) \geq n$.

Since $\bar I$ is nilpotent we have that $K_2^{(n)}=\{1\}$ for all large $n$.
Clearly we have $K_2^{(0)}=K_2(\bar{R})$, then
$$
K_2(\bar R)=K_2^{(0)}\geq K_2^{(1)} \geq \cdots \geq \{1\}
$$
is a filtration of $K_2(\bar{R})$ terminating at $\{1\}$.

\begin{lemma} \label{simplify}
Let $v$ be a unit of $\bar R$, $k \in \N$ and $p,q,r \in \bar I$.

(i1)\  If $l(p)+l(q) \geq n$ then
$$
\{v, 1-pqv\}\{1-qv,1-p\} \{1-pv,1-q\} \equiv 1 \mb{ mod } K_2^{(n+1)},
$$

(i2)\  If $l(p)+l(q)+l(r) \geq n$ then
$$
\{1-p,1-qr\} \{1-q,1-pr\}\{1-r,1-pq\} \equiv 1 \mb{ mod }  K_2^{(n+1)}.
$$
\end{lemma}
\textbf{Proof:}
Using the basic properties of the Steinberg symbols we have
$$
\begin{array}{l@{}c@{}l}
\displaystyle
\left\{  -\frac{1-qv}{1-p},  \frac{1-pqv}{1-p} \right\}
& {}={}&
\displaystyle
\left\{  {1-qv},  {1-pqv} \right\}
\left\{  {1-qv},  \frac{1}{1-p} \right\} \times \\
& & \quad \quad \times
\displaystyle
\rule[18pt]{0pt}{0pt}
\left\{  -\frac{1}{1-p},  {1-pqv} \right\}
\left\{  -\frac{1}{1-p},  \frac{1}{1-p} \right\} = \\
& {}={} &
\displaystyle
\rule[15pt]{0pt}{0pt}
\left\{  {1-qv}, {1-pqv} \right\}
\left\{  {1-p} , {1-qv}  \right\} \times \\
& & \quad \quad \times
\displaystyle
\rule[15pt]{0pt}{0pt}
\left\{  {1-p}, {1-pqv} \right\}^{-1}\left\{-1, {1-pqv} \right\}.
\end{array}
$$
The first and the third multiplier are in $K_2^{(n+1)}$. The forth multiplier is the identity if the characteristic of $\bar R/\bar I$ is $2$ 
and in $K_2^{(n+1)}$ otherwise.
This gives us
$$
\left\{  -\frac{1-qv}{1-p},  \frac{1-pqv}{1-p} \right\} \equiv
\left\{  {1-p} , {1-qv} \right\} \mb{ mod } K_2^{(n+1)}
$$
Applying this reduction in the identity (\ref{I1}) leads to the (i1).
In the same way we can obtain (i2) from the identity (\ref{I2}).
$\square$

\medskip

Identity (i2) from Lemma \ref{simplify} has the following:
\begin{corollary}\label{product}
Let $r,q_1,q_2,...,q_t \in \bar I$ and suppose that
$l(r)+\sum_i l(q_i) \geq n$.
Then we can write $\{1-q_1q_2\cdots q_t,1-r\}$ as
$$
\{1-q_1q_2\cdots q_t,1-r\}=\prod_{i=1}^t \{1-q_i,1-u_i\} \mb{ mod }  K_2^{(n+1)}
$$
for some
$u_i\in \bar I$ with $l(q_i)+l(u_i)\geq n$ for each $i=1,2,\dots,t$.
\end{corollary}

Now, fix $\theta$ to be an element in $\bar R$ such that the image of $\theta$ generates
the multiplicative group
$\bar R/\bar I$, and let $a_j, \ j=1,2,\dots,k+1$  be some fixed set of generators of the ideal $\bar I$.
\medskip

We \textbf{claim} that the identities above and the fact that the ideal $\bar I$ is
nilpotent allow us to write (not uniquely) every
Steinberg symbol $\{a,b\}$ as a product
$$
\{a,b\} = \{\theta, * \}\prod_{j=1}^{k+1}\{1-a_j,*\}.
$$

Let $T$ denote the following subset of $K_2(\bar R)$:
$$
T:=\left\{\{\theta, s_0\}\prod_{j=1}^{k+1}\{1-a_j, s_j\} \bigg| \ s_j \in \bar R^*, s=0,1,\dots,k \right\}.
$$
Bimultiplicativity of $\{-,-\}$ implies that $T$ is a subgroup of $K_2(\bar R)$.
\medskip

We shall prove by induction on $n$ that any symbol $\{a,b\}$
can be written as $\{a,b\}=tu$ with $t \in T$ and
$u\in K_2^{(n)}$. As $K_2^{(n)}=\{1\}$ for large
enough $n$ this will prove the above claim and Theorem \ref{k2}.

Now, every unit $a\in \bar R$ is congruent to
$\theta^l$ for some integer $l$, and thus $a=\theta^l a'$ with $a'\in 1+I$. Then
$$
\{a,b\}=\{\theta^l,b\} \{a',b\} = \{\theta,b^l\} \{a',b\}
$$
and $\{a',b\} \in K_2^{(1)}$, giving the first step $n=1$ of the induction.
\medskip

Suppose our claim has been proved for $n-1$. We need to express a symbol
$\{a,b\}$ with $l(1-a)+l(1-b)=n$ 
in the required form modulo $K_2^{(n+1)}$.
The same argument from the base step $n=1$ gives that it is sufficient
to consider the case $a=1-r, b=1-q$ with $r,q \in I$ and $l(r)+l(q)=n$.
The ideal $\bar I$ is 
generated by $a_i$ modulo $\bar I^2$. Therefore $r$ may be written as
$$
r\equiv  \sum_{i=1}^s v_i a_{m_{1,i}}a_{m_{2,i}}\cdots a_{m_{t,i}} 
\quad \textrm{mod } \bar I^{t+1},
$$
where $t=l(q)$ and the $v_i$ are units in $\bar R$.
Now
$$
1-r \equiv \prod_{i=1}^s ( 1-v_i a_{m_{1,i}}a_{m_{2,i}}\cdots a_{m_{t,i}})
\quad \textrm{mod } \bar I^{t+1}
$$
giving that
$$
\{a,b\} \equiv \prod_{i=1}^s \{ 1-v_i a_{m_{1,i}}a_{m_{2,i}}\cdots a_{m_{t,i}},b \}
\quad \textrm{mod }\ K_2^{(n+1)} .
$$

By Corollary \ref{product} each term
$\{1-v_i a_{m_{1,i}}a_{m_{2,i}} \cdots a_{m_{t,i}},b \}$
in the product on the right hand side is equal to
$$
\{1-v_ia_{m_{1,i}},1-u_{i,1}\}\prod_{j=2}^t\{ 1-a_{m_{j,i}},1-u_{i,j}\}
 \quad \mathrm{ mod } \ K_2^{(n+1)}
$$
for the appropriate $u_{i,j} \in \bar I$, $(j=1,2,\dots,t)$.
These symbols are all in the required from, except maybe the first one:
$\{1-v_ia_{m_{1,i}},1-u_{i,1}\}$.
\medskip

Recall that by Corollary \ref{product} we have
$l (a_{m_{1,i}})+l(u_{i,1}) \geq n$
and therefore (i1) of Lemma \ref{simplify}
(with $v=v_i, p=u_{i,1}, q=a_{m_{i,1}}$ ) gives
$$
\{1-v_ia_{m_{1,i}}, 1-u_{i,1}\}
\equiv
\{v_i,1-u_{i,1}a_{m_{i,1}}\}^{-1}
\{1-a_{m_{i,1}},1-v_iu_{i,1}\}
\quad \mathrm{ mod }\ K_2^{(n+1)}.
$$
The unit $v_i$ can be written as
$v_i \equiv \theta^{-w_i}$ mod $\bar I$ for some $w_i\in \Z$.
Now
$$
\{v_i,1-u_{i,1}a_{m_{i,1}}\}^{-1}\equiv \{\theta, 1-u_{i,1}a_{m_{i,1}}\}^{w_i}
\equiv \{\theta, (1-u_{i,1}a_{m_{i,1}})^{w_i}\}
\quad \mathrm{ mod }\ K_2^{(n+1)}
$$
and we have expressed
$\{a,b\}$ as an element of $T$ modulo $K_2^{(n+1)}$.
\bigskip

This completes the induction and proves our claim and Theorem \ref{k2}. $\square$

\section{On the profinite completion of $G$}
\label{CSP}
In this section we shall prove Theorem \ref{BpEG} and several
basic properties of the congruence kernel of $G$.
\medskip

As noted in Section \ref{outline} Theorem \ref{k2} together with
Lemma \ref{congimage} imply that any finite image $G/H$ of $G$ is
a product of at most $N:=(3d-d-2)/2 + 18(k+2)$ of the images
$E_{i,j}H$ of the elementary groups $E_{i,j}$ in some fixed order. We make this
statement more precise:

There is a map
$s \mapsto E(s)=E_{i_s,j_s} \in \{E_{i,j} |\ i\not =j \}$ for $s=1,2,\dots,N$ such that
$$
G/H =\left(\prod_{s=1}^N E(s) \right) H
$$
for any normal subgroup $H$ of finite index in $G$.
\medskip

Recall that $\widehat G$ is the profinite completion of $G$ and
$\widehat E_{i,j}\simeq \widehat R$ is the closure of $E_{i,j}$ in $\widehat G$
(where $\widehat{R}$ is the profinite completion of $R$).

Set $\widehat E(s)= \widehat E_{i_s,j_s} \subseteq \widehat G$.
It is a closed subgroup of $\widehat G$ and the identity above shows that
$\Pi:=\prod_{s=1}^N \widehat E(s)$ is dense in $\widehat{G}$.
However the set $\Pi$ is a finite product of closed sets in a compact Hausdorff group
and is therefore closed. Hence $\Pi=\widehat G$ and Theorem \ref{BpEG} is proved. $\square$
\bigskip

Let $\widehat G_c=\widehat{G_{d,k}}_c$ be the congruence completion of $G=G_{d,k}$, i.e.,
$$
\widehat G_c=
\lim_{\substack{ \longleftarrow \\ U \triangleleft R}}G/G(U)\simeq
\lim_{ \substack{ \longleftarrow \\ U\triangleleft R}}
G(R/U)\simeq G(\widehat R).
$$
We have natural surjection $p:\ \widehat G \rightarrow \widehat G_c$.

Define $\mathcal{C}_{d,k}:=\ker p$, this is called the \emph{congruence kernel} of $G_{d,k}$.
From the argument in Section \ref{outline} we see that
$$
\mr{C}_{d,k} \simeq  \lim_{\substack{ \longleftarrow \\ U \triangleleft R}}G(U)/EL(d,U)=
\lim_{\substack{ \longleftarrow \\ U \triangleleft R}}SK_1(R,U;d).
$$

\medskip

The centrality of $K_2(R/U;d)$ in $\mathrm{St}_d(R/U)$ (see~\cite{DS2}) gives that
$\mr{C}_{d,k}$ is central in $\widehat{G}$.
We thus have the first part of the following Theorem:
\begin{theorem}
\label{C}
Suppose $d\geq 3$ and $k\geq 1$.
The congruence kernel $\mr{C}_{d,k}$ is a central closed subgroup of
$\widehat G$. Moreover for $d_1,d_2>k+2$ the abelian groups $\mr{C}_{d_1,k}\simeq
\mr{C}_{d_2,k}$ are isomorphic and not finitely generated as profinite groups.
\end{theorem}

To prove the second claim note that the polynomial ring $R$ satisfies the Bass
stable range condition $\text{sr} (R)\leq k+2$ (see Definition 3.1 and Theorem 3.5 of Chapter V in \cite{bass}).

A result of Vaserstein \cite{V}, improving on \cite{bass} gives an isomorphism
\begin{equation} \label{stabilization}
SK_1(R,U;d_1)\simeq SK_1(R,U;d_2) \simeq SK_1(R,U)
\end{equation}
for $d_i>k+2$.
\medskip
Recall the exact sequence from Section \ref{symbols}:
\[ K_2(R) \rightarrow K_2(R/U) \rightarrow SK_1(R,U) \rightarrow 1 \]
A classical result by Quillen \cite{quillen} is that $K_2(R)\simeq K_2(\mathbb{Z}) \simeq \{\pm 1\}$.
On the other hand there are various examples of finite images $R/U$ such that $K_2(R/U)$ has large rank,
e.g.\ if $R=\mathbb{Z}[x]$ and $U=\langle p^2, x^{p^l}\rangle$ then $K_2(R/U)$
is an elementary abelian $p$-group of rank $\geq l$
by \cite{DS2}, proof of Theorem 2.8.
We conclude that $SK_1(R,U)$ can have arbitrary large rank and by (\ref{stabilization}) the same is true for the
finite images $SK_1(R,U;d)$ of $\mr{C}_{d,k}$. Therefore the congruence kernels $\mr{C}_{d,k}$
are not finitely generated for $d>k+2$ and $k\geq 1$, which proves
the second claim. 

\medskip

We remark that Theorem \ref{k2} can be restated in the following form:
\begin{remark}
The congruence kernel $\mathcal{C}_{d,k}$ has finite width with respect to the Steinberg symbols
$\{a,b\}\in \mathrm{St}_d(\widehat R)$
with arguments $a,b \in \widehat{R}^*$.
\end{remark}

Theorem~\ref{C} implies the following corollary which gives a negative answer to
Question 1.27 from~\cite{LubZuk}:
\begin{corollary}
Provided $d > k+2$ and
$k\geq 1$ the congruence completion $\widehat{G_{d,k}}_c$ of the universal lattice $G_{d,k}$
is a profinite group with property $\tau$ which is not finitely presented.
\end{corollary}


\noindent
\texttt{Martin Kassabov, \\ Cornell University, Ithaca, NY 14853-4201, USA. \\
\emph{e-mail:} kassabov@math.cornell.edu}

\medskip

\noindent
\texttt{Nikolay Nikolov, \\ New College, Oxford, OX1 3BN, UK. \\
\emph{e-mail:} nikolay.nikolov@new.ox.ac.uk}


\begin{thebibliography}{99}

\bibitem{bass} H. Bass, \emph{Algebraic K-theory},
W.A.Benjamin, Inc, New York-Amsterdam, 1968. (Chapter V)


\bibitem{Bur} M. Burger,
Kazhdan constants for ${\rm SL}_3({\rm Z})$,
\emph{J. Reine Angew. Math., 413 (1991), 36--67.}

\bibitem{CK} D. Carter, G. Keller,
Bounded elementary generation of ${\rm SL}\sb{n}({\cal O})$.
\emph{Amer. J. Math 105 (1983), no.3,  673--687.}

\bibitem{DS1}  R. K. Dennis, M. R. Stein,
$K\sb{2}$ of discrete valuation rings.
\emph{ Advances in Math. 18 (1975), no. 2, 182--238.}

\bibitem{DS2} R. K. Dennis, M. R. Stein,
$K\sb{2}$ of radical ideals and semi-local rings revisited.
\emph{Algebraic $K$-theory, II: "Classical" algebraic $K$-theory and connections with arithmetic}
(Proc. Conf., Battelle Memorial Inst., Seattle, Wash., 1972),
Lecture Notes in Math. Vol. 342, Springer, Berlin, 1973, pp 281-303.

\bibitem{DV}R. K. Dennis, L. N. Vaserstein, On a question of M. Newman on the number of commutators.
\emph{J. Algebra 118 (1988), no. 1, 150--161.}

\bibitem{gromov} M. Gromov,
Random walks in random groups,
\emph{GAFA, Geom. funct. anal. 13:1 (2003), 73--148.}

\bibitem{K} M. Kassabov,
Kazhdan constants for $\rm{SL}_d(\mathbb{Z})$,
\texttt{arXiv:math.GR/0311487}, \emph{to appear in Internat. J. Algebra and Comput.}

\bibitem{KSL3k} M. Kassabov,
Universal lattices and unbounded rank expanders,
\texttt{arXiv:math.GR/0502237}.

\bibitem{kazhdan} D. A. Kazhdan,
On the connection of the dual space of a group with the structure of its closed subgroups,
\emph{Funkcional.
Anal. i Priloz. 1 (1967), 71-74.}

\bibitem{expanderbook} A. Lubotzky,
\emph{Discrete groups, expanding graphs and invariant measures},
Progr. Math. 125, Birkh\"{a}user, Boston, 1994.

\bibitem{LubZuk} A. Lubotzky, A. $\hat{\rm{Z}}$uk,
\emph{On Property ($\tau$)}, book.

\bibitem{milnor} J. Milnor,
\emph{Introduction to Algebraic K-Theory}, Princeton University Press, 1971.

\bibitem{quillen} D. Quillen, Higher algebraic $K$-theory, \emph{Algebraic $K$-theory, I:
Higher $K$-theories}
(Proc. Conf., Battelle Memorial Inst., Seattle, Wash., 1972), Lecture Notes in Math., Vol. 341,
Springer, Berlin 1973, pp. 85-147.

\bibitem{ZS} P. Samuel, O. Zariski,
\emph{Commutative Algebra II}, Graduate Texts in Mathematics, Vol. 29, Springer-Verlag, 1975.

\bibitem{YSh} Y. Shalom,
Bounded generation and Kazhdan's property (T),
\emph{Publ. Math. IHES, 90 (1999), 145-168.}

\bibitem{stein} M. R. Stein, Surjective stability in dimension $0$ for $K\sb{2}$ and related functors.
\emph{Trans. Amer. Math. Soc. 178 (1973), 165-191.}

\bibitem{suslin} A. A. Suslin,
The structure of the special linear group over rings of polynomials. (Russian),
\emph{ Izv. Akad. Nauk SSSR Ser. Mat. 41 (1977), no. 2, 235--252, 477.} MR0472792

\bibitem{V} L. N. Vaserstein,
On the stabilization of the general linear group over a ring.
\emph{Mat. Sb. (N.S.) 79 (121) 405--424 (Russian); translated as Math. USSR-Sb 8 1969 383--400.}

\bibitem{zuk} A. $\hat{\mathrm{Z}}$uk, Property T and Kazhdan constants for discrete groups,
\emph{Geom. Funct. Anal. 13,  no.3,
643-670}

\end{thebibliography}
\end{document}